\input amstex
\documentstyle{amsppt}
\NoBlackBoxes

\TagsOnRight

\def\cal{\Cal}
\def\AA{{\cal A}}

\def\HH{{\cal H}}
\def\MM{{\cal M}}

\def\SS{{\cal S}}

\def\LL{{\cal L}}

\def\PP{{\cal P}}
\def\Z{{\Bbb Z}}
\def\C{{\Bbb C}}
\def\R{{\Bbb R}}

\def\Q{{\Bbb Q}}
\def\e{{\epsilon}}

\def\n{\noindent}
\def\part{{\partial}}
\def\dudtau{{\part u \over \part \tau}}

\rightheadtext{Length minimizing property} \leftheadtext{
Yong-Geun Oh }

\topmatter
\title
Length minimizing property, Conley-Zehnder index and
$C^1$-perturbations of Hamiltonian functions
\endtitle
\author
Yong-Geun Oh\footnote{Partially supported by the NSF Grant \#
DMS-0203593,  Vilas Research Award of University of Wisconsin and
by a grant of the Korean Young Scientist Prize \hskip8.5cm\hfill}
\endauthor
\address
Department of Mathematics, University of Wisconsin, Madison, WI
53706, ~USA \& Korea Institute for Advanced Study, 207-43
Cheongryangri-dong Dongdaemun-gu Seoul 130-012, KOREA
\endaddress

\abstract The main purpose of this paper is to study the length
minimizing property of Hamiltonian paths on closed symplectic
manifolds $(M,\omega)$ such that there are no spherical homology
class $A \in H_2(M)$ with
$$
\omega(A) > 0 \quad \text{and} \quad  -n \leq c_1(A) < 0,
$$
which we call {\it very strongly semi-positive}. We introduce the
notion of {\it positively $\mu$-undertwisted} Hamiltonian paths
and prove that any positively undertwisted quasi-autonomous
Hamiltonian path is length minimizing in its homotopy class as
long as it has a fixed maximum and a fixed minimum point that are
generically under-twisted. This class of Hamiltonian can have
non-constant large periodic orbits. The proof uses the chain level
Floer theory, spectral invariants of Hamiltonian diffeomorphisms
and the argument involving the thick and thin decomposition of
Floer's moduli space of perturbed Cauchy-Riemann equation. And
then based on this theorem and some closedness of length
minimizing property, we relate the Minimality Conjecture on the
very strongly semi-positive symplectic manifolds to a
$C^1$-perturbation problem of Hamiltonian functions on general
symplectic manifolds, which we also formulate here.
\endabstract


\endtopmatter

\document

\bigskip

\head \bf Contents \endhead

\n 1. Introduction
\smallskip
\n 2. The Hofer topology and the length minimizing property
\smallskip
\n 3. Chain level Floer theory and spectral invariants
\smallskip
\n 4. Canonical fundamental Floer cycles
\smallskip
\n 5. Proof of Theorem $\text{A}'$ and B
\smallskip
\n 6. Minimality and the $C^1$ Perturbation Conjecture

\bigskip

\bigskip

\head  \bf \S1. Introduction
\endhead
The celebrated Hofer's norm of Hamiltonian diffeomorphisms
introduced in [H] is defined by
$$
\|\phi\| = \inf_{H\mapsto \phi} \|H\| \tag 1.1
$$
where $H \mapsto \phi$ means that $\phi= \phi_H^1$ is the time-one
map of Hamilton's equation
$$
\dot x = X_H(x)
$$
and the norm $\|H\|$ is  defined by
$$
\|H\| = \int _0^1 \text{osc }H_t \, dt = \int_0^1(\max H_t -\min
H_t)\, dt. \tag 1.2
$$
This induces a distance function, the so called the Hofer
distance, on $Ham(M,\omega)$ which in turn defines a topology on
$Ham(M,\omega)$. We call this topology the {\it Hofer topology} on
$Ham(M,\omega)$.

Our convention of the definition of Hamiltonian vector field will
be
$$
X_h \rfloor \omega = dh \tag 1.3
$$
for a smooth function $h$ on $M$. From now on we will always
assume $M$ is closed and, that the Hamiltonian functions are
normalized, unless otherwise said, so that
$$
\int_M H_t\, d\mu = 0
$$
where $d\mu$ is the Liouville measure of $(M,\omega)$.

The main purpose of the present paper is to study the length
minimizing property of Hamiltonian paths without the standard
condition that the corresponding Hamiltonian has {\it no}
non-constant periodic orbits, but instead with some {\it
topological} condition on the Conley-Zehnder indices of the
periodic orbits. Our motivation for this attempt is to study the
following prominent question posed by Polterovich [Conjecture
12.6.D, Po] and by Lalonde-McDuff-Slimowitz [LM2], [MS] in this
regard.

\proclaim{[Minimality Conjecture]} Any autonomous Hamiltonian path
that has no contractible periodic orbits of period less than equal
to one is Hofer-length minimizing in its homotopy class with fixed
ends.
\endproclaim

Our attempt has led us to the following definition.

\definition{Definition 1.1} We call a symplectic manifold
$(M,\omega)$ {\it very strongly semi-positive} if it does not
carry any spherical homology class $A \in H_2(M,\Z)$ such that
$$
\omega(A) > 0, \quad -n \leq c_1(A) < 0. \tag 1.4
$$
We denote by
$$
\Sigma : = \min\{|c_1(A)| \mid c_1(A) \neq 0 \}
$$
and call the {\it minimal Chern number} of the symplectic
manifold.
\enddefinition
This class of symplectic manifolds include all weakly exact
symplectic manifolds, all (positively) monotone symplectic
manifolds, and negatively monotone symplectic manifolds with
$\Sigma \geq n +1 = \frac{1}{2}\dim M + 1$. In particular, it
includes any Fano manifold, e.g., $\C P^n$.

We would like to emphasize that although the very strongly
semi-positive condition looks very much alike and slightly
stronger than the usual semi-positivity condition [HS] in which
the condition on $c_1$ in (1.4) is replaced by
$$
-n + 3 \leq c_1(A) < 0
$$
in relation to the transversality question of pseudo-holomorphic
curves, the origin of our very strongly semi-positive condition
here is different from that of the semi-positivity condition. In
particular, the enhanced machinery of virtual moduli cycles is
irrelevant to our requirement and will {\it not} help removing
this condition from the statement of the main result below.
However under the very strongly semipositivity assumption, the
technical aspect of the Floer theory in this paper is `elementary'
in that it does not require such enhanced machinery. Here we avoid
using the more natural name like `strongly semi-positive' because
it has been already used in [En] for the case where the $c_1$
condition in (1.4) is replaced by $-n + 2 \leq c_1(A) < 0$.

We now recall some basic definitions in the study of Hamiltonian
diffeomorphisms and their paths. Two Hamiltonians $G$ and $F$ are
called equivalent if there exists a smooth family $\{F^s\}_{0\leq
s\leq 1}$ with $F^0 = G, \, F^1 = F$ such that
$$
\phi^1_{F^s} = \phi^1_G
$$
for all $s \in [0,1]$. We denote $G \sim F$ in that case and say
that two Hamiltonian paths $\phi^t_G$ and $\phi^t_F$ are homotopic
to each other with fixed ends, or just homotopic to each other
when there is no danger of confusion.

\definition{Definition 1.2} A Hamiltonian $H$ is called {\it
quasi-autonomous} if there exists two points $x_-, \, x_+ \in M$
such that
$$
H(x_-,t) = \min_x H(x,t), \quad H(x_+,t) = \max_x H(x,t)
$$
for all $t\in [0,1]$.
\enddefinition

To state the result on the length minimizing property in this
paper and for the future purpose, the following general definition
seems to be useful. This is a {\it topological} undertwistedness
while the ones used in [En, MS] is {\it dynamical}. We denote by
$[z,w]$ an element of the standard $\Gamma$-covering space
$\widetilde\Omega_0(M,\omega)$ [HS] of the contractible loop free
loop space $\Omega_0(M)$, and by $\mu_H([z,w])$ its Conley-Zehnder
index.

\definition{Definition 1.3} Let $H$ be any Hamiltonian
and $z$ be a contractible one-periodic orbit of its Hamiltonian vector field
$X_H$.
\roster
\item
We say that $z$ is {\it positively $\mu_H$-undertwisted} if it allows a
bounding disc $w_z$ such that
$$
-n \leq \mu_{H}([z,w_z]) \leq n \quad\text{and}\quad \int
w_z^*\omega \geq 0
$$
for $[z,w_z]$. If this holds for all contractible one-periodic
orbits of $H$, then we call $H$ {\it positively
$\mu$-undertwisted}.
\item
If $\int w_z^*\omega \geq - \e$ holds for $\e >0$ instead, then we
say the periodic orbit $z$ {\it $\e$-positively $\mu_H$-undertwisted},
and $H$ {\it $\e$-positively $\mu$-undertwisted} respectively.
\endroster
\enddefinition

Examples of positively $\mu$-undertwisted Hamiltonians will be any
{\it slow} autonomous Hamiltonian in the sense of [En, MS], and
any $C^2$-small perturbations of nondegenerate slow autonomous
Hamiltonians are $\e$-positively $\mu$-undertwisted for some
sufficiently small $\e
> 0$.

\proclaim{Theorem A} Let $(M,\omega)$ be very strongly
semi-positive. Suppose $H$ is a quasi-autonomous Hamiltonian such
that \roster
\item it has a fixed global maximum and a fixed global minimum
point that are generically under-twisted.
\item $H$ is positively $\mu$-undertwisted and nondegenerate in
the Floer theoretic sense
\endroster
Then the Hamiltonian path $\phi_H^t$ is length minimizing in its
homotopy class with fixed ends.
\endproclaim
We would like to emphasize that the Hamiltonian $H$ can have
`large' periodic orbits as long as they are positively
$\mu_H$-undertwisted.

In fact, the proof of this theorem shows that the condition (2)
can be weakened to the following

\proclaim{Theorem $\text{A}'$} Let $(M,\omega)$ be very strongly
semi-positive. Let $H$ satisfy just (1) in Theorem A and $e_H >0$
be a constant depending only on `the local behavior' of $H$ at the
maximum and minimum points introduced in Lemma 5.4 later. Suppose
that there exists a sequence of $\e \to 0$ such that $H$ allows a
$C^1$-small perturbation $H'$ that is quasi-autonomous and
satisfies the condition that
\smallskip

$(2^\prime)$ $H'$ is $\e$-positively $\mu$-undertwisted and
nondegenerate in the Floer theoretic sense.
\smallskip

\n Then the Hamiltonian path $\phi_{H}^t$ is length minimizing in
its homotopy class with fixed ends.
\endproclaim
This perturbation result is the precise formulation of the remark
mentioned right after [Theorem I, Oh3]. We refer to section 5 for
the precise meaning of `the local behavior' in the statement.

We would like to emphasize that the perturbation in Theorem
$\text{A}'$ is assumed to be $C^1$-small, not necessarily
$C^2$-small. In particular, there is no a priori relation between
the Conley-Zehnder indices of the nearby periodic orbits of $H$
and $H'$ at all. However combined with the closedness of the
length minimizing property (Theorem 3.1) of Hamiltonian paths
under the Hofer topology, this will be an important point in
application to the Minimality Conjecture which lead us to the
$C^1$ Perturbation Conjecture later.

Our proof of Theorem A and $\text{A}'$ will be based on the chain
level Floer theory from [Oh3,5] using the scheme developed in
[Oh6] via the usage of the spectral invariant $\rho(H;1)$. Many
arguments are combination of those from [Oh3,6].

First to recall this criterion, we rewrite the Hofer norm into
$$
\|H\| = E^-(H) + E^+(H)
$$
where $E^\pm$ are the negative and positive parts of the Hofer
norms defined by
$$
\align
E^-(H) & = \int_0^1 -\min H\, dt \\
E^+(H) & = \int_0^1 \max H\, dt.
\endalign
$$
These are called the {\it negative Hofer-length} and the {\it
positive Hofer-length} of $H$ respectively. We note
$$
E^+(H) = E^-(\overline H)
$$
where $\overline H$ is the Hamiltonian generating
$(\phi_H^t)^{-1}$ defined by
$$
\overline H(t,x) = - H(t, \phi_H^t(x)).
$$
Therefore we will focus only on the semi-norm $E^-$. According to
the criterion [Theorem III, Oh6], Theorem $\text{A}'$ will be an
immediate corollary of Theorem 3.1 [Lemma 5.1, Oh3] and the
following theorem in whose proof the very strongly semi-positive
condition of $(M,\omega)$ enters in an essential way. We refer to
section 3 for a brief outline of the construction from [Oh5] of
the spectral invariant $\rho(H;1)$ on non-exact symplectic
manifolds.

\proclaim{Theorem B} Let $(M,\omega)$ and $H'$ be as in Theorem
$\text{A}'$. Then we have
$$
\rho(H';1) = E^-(H').
$$
The same holds for the inverse Hamiltonian $\overline H'$.
\endproclaim

The main new ingredient in the proof of Theorem B besides the
scheme used in [Oh3,6] is our usage of the argument in section 6
that is based on the ``thick and thin'' decomposition of the Floer
moduli space. Similar argument was previously used by the author
[Oh1] in the context of Lagrangian intersection Floer homology
theory for its application to the Maslov class and construction of
the corresponding spectral sequence.

We recall that the Minimality Conjecture was answered
affirmatively for the weakly exact case or of the case of surfaces
by Lalonde-McDuff much earlier [LM2] by a different method. Our
proof also works for the surface case without any further
requirement on $G$ other than those in the conjecture, but does
require (1) even for the weakly exact case when $\dim M \geq 4$.
We, however, refer to section 2 and 6 for some discussion on how
one might be able to improve this point in the general context. It
would be very interesting to see if the very strongly
semi-positive condition is an essential condition or just a
technical artifact of our Floer theoretic approach. Roughly
speaking, the very strongly semi-positive condition rules out a
possible `quantum contribution' to the minimization process. It
seems to be of fundamental importance to understand how the
quantum contribution affects the length minimization process in
general. We hope to come back to this issue in the future.

On the other hand, based on Theorem $\text{A}'$ and the following
$C^1$ Perturbation Conjecture on the Hamiltonians and Theorem 3.1
[Lemma 5.1, Oh3], we believe that for the autonomous Hamiltonian
$G$ appearing in the Minimality Conjecture the condition (2) can
be replaced by the condition ``$G$ has no non-constant periodic
orbit'' at least for the very strongly semi-positive case.

\medskip

\proclaim{[$C^1$ Perturbation Conjecture]} Let $(M,\omega)$ be any
symplectic manifold. Suppose that $G$ is an autonomous Hamiltonian
with a global maximum and a minimum that are generically
under-twisted, and has no non-constant periodic orbits of period
less than equal to one. Let $\e_G > 0$ be the positive constant in
Theorem $\text{A}'$. Then there is a $C^1$-small perturbation $H'$
of $G$ such that $H'$ is $\e$-positively $\mu$-undertwisted where
$\e$ can be made arbitrarily small depending only on $G$.
\endproclaim

We refer to section 6 for a more precise formulation of the
conjecture and for some further discussion. If this conjecture is
true, then the Minimality Conjecture will follow from Theorem
$\text{A}'$ at least for the very strongly semi-positive case.
Based on the way how the proof goes in section 5, we suspect that
the conjecture is not true in general.

The organization of the paper is in order. Section 2 through 4
deal with the cases of general symplectic manifolds. The very
strongly semi-positivity condition enters only in section 5. In
section 2, we recall a theorem proven in [Oh7] which states that
the length minimizing property is {\it closed} under the { \it
Hofer topology} of Hamiltonian paths. This theorem is a stronger
statement than the one stated in [Lemma 5.1,Oh3].

In section 3, we recall the basic chain level Floer theory, and
briefly outline the construction from [Oh3,5] of spectral
invariants, especially the one $\rho(H;1)$ on {\it non-exact}
symplectic manifolds, and state their basic properties. See also
[Oh2], [Sc] for the construction in the exact case.

In section 4, we recall the notion of {\it the canonical
fundamental Floer cycle} that was introduced in [Oh3,6] and
compute its {\it level}. Since the scheme of the proof in
[Proposition 4.3, Oh6] will play an essential role in the proof of
Theorem A, for the reader's convenience, we duplicate the proof
here.

In section 5, we restrict to the very strongly semi-positive case
and give the proof of Theorem $\text{A}'$ and B by proving that
the above canonical fundamental cycle is {\it tight} in the sense
of Definition 5.2 [Definition 4.2, Oh6] under the condition of
Theorem $\text{A}'$ and B which will in turn prove Theorem B and
hence Theorem $\text{A}'$.

In section 6, we formulate a precise version of the above $C^1$
Perturbation Conjecture and explain how the conjecture together
with Theorem $\text{A}'$ would imply the Minimality Conjecture for
the very strongly semi-positive $(M,\omega)$.

In the appendix, we provide complete details of a proof of the
index formula
$$
\mu_H([z,w]) = \mu_H([z,w']) + 2c_1(w\#\overline w') \tag 1.5
$$
where we emphasize the sign `$+$' in front of the Chern number
term in this formula. An incorrect formula with the negative sign
was written in [section 6.1, En] and [section 7, Oh3]. We also
clarify the correct formula under the other commonly used package
of conventions as in [Po] and others. Having a correct sign did
not play any significant role in the works in [En], [Oh3] or in
other previous literature on the symplectic Floer homology theory.
Since our definition of  the {\it positively
$\mu$-undertwistedness} and the proof of the main theorems in this
paper crucially depend on having the correct sign in this formula
and we cannot locate any reference that contains a complete proof
of the formula in any convention, we give complete details of the
proof essentially from the scratch based on the definition of the
Conley-Zehnder index given in [CZ], [SZ] on $\R^{2n}$ for the reader's
convenience. The reference [HS] contains the formula (1.5) but
their conventions do not completely agree with ours and we do not
feel safe just to quote the formula, especially when we are not
completely sure of what conventions the authors of [HS] are using.
There have already been more than one instances of an incorrect
formula written in the literature as in [En] and [Oh3].

We thank the unknown referee of the originally submitted version
of the paper [Oh6] for pointing out that the Minimality Conjecture
does not follow from the main result in the original version of
[Oh6]. This remark prompted us to carefully look over the delicate
points of appearance of ``small'' periodic orbits, which we had
overlooked in [Oh3] (see Erratum to [Oh3] for clarification of
these points). The final writing has been carried out while we are
visiting Korea Institute for Advanced Study during the winter of
2003-2004. We thank KIAS for its financial support and excellent
research atmosphere during our stay.

\head{\bf \S 2. The Hofer topology and the length minimizing
property}
\endhead

In this section, we first recall the notion of {\it Hofer
topology} on the space of Hamiltonian paths, and the cloosedness
of the length minimizing property of Hamiltonian paths under the
topology proven in [Oh3,7]. Our presentation of the Hofer topology
closely resembles that of the {\it Hamiltonian topology} that we
introduce in [Oh7], which is however equivalent to the usual
description of the Hofer topology in the literature.

We first recall the definition of Hamiltonian diffeomorphisms and
smooth Hamiltonian paths.

\definition{Definition 2.1}
(i) A $C^\infty$ diffeomorphism $\phi$ of $(M,\omega)$ is a {\it
Hamiltonian diffeomorphism} if $\phi = \phi_H^1$ for a $C^\infty$
function $F: \R \times M \to \R$ such that
$$
F(t+1,x) = F(x) \tag 2.1
$$
for all $(t,x) \in \R \times M$. Having this periodicity in mind,
we will always consider $F$ as a function on $[0,1] \times M$.
\smallskip

\n (ii) A {\it Hamiltonian path} $\lambda: [0,1] \to
Ham(M,\omega)$ is a smooth map
$$
\Lambda: [0,1] \times M \to M
$$
such that \roster
\item its derivative $\dot\lambda(t) = \frac{\part \lambda}{\part t}\circ
(\lambda(t))^{-1}$ is Hamiltonian, i.e., the one form
$\dot\lambda(t) \, \rfloor \omega$ is exact for all $t \in [0,1]$.
We call the normalized function $H: \R \times M \to \R$ the {\it
generating Hamiltonian} of $\lambda$ if it satisfies
$$
\lambda(t) = \phi_H^t(\lambda(0)) \quad \text{or equivalently
$dH_t = \dot\lambda(t) \, \rfloor \omega$}.
$$
\item
$\lambda(0):=\Lambda(0, \cdot): M \to M$ is a Hamiltonian
diffeomorphism, and so is for all $\lambda(t) = \Lambda(t, \cdot),
\, t\in \R$.
\endroster
We denote by $\PP(Ham(M,\omega))$ the set of Hamiltonian paths
$\lambda: [0,1] \to Ham(M,\omega)$, and by $\PP(Ham(M,\omega),id)$
the set of  $\lambda$ with $\lambda(0) = id$.
\enddefinition

Here $\phi_H^1$ is the time-one map of the Hamilton equation
$$
\dot x = X_H(x). \tag 2.2
$$
We will always denote by $\phi_H$ the corresponding Hamiltonian
path
$$
\phi_H: t \mapsto \phi_H^t
$$
starting from the identity, and by $H \mapsto \phi$ when $\phi =
\phi_H^1$. In the latter case, we say that the diffeomorphism
$\phi$ is generated by the Hamiltonian $H$.

The {\it Hofer length} of the Hamiltonian path $\lambda \in
\PP(Ham(M,\omega))$ with $\lambda(t) = \phi_H^t(\lambda(0))$ is
given by the norm of its generating Hamiltonian $H$ defined by
$$
\text{leng}(\lambda) = \|H\| = \int_0^1(\max_x H_t -\min_x H_t)\,
dt.
$$
We also denote by
$$
ev_1:\PP(Ham(M,\omega),id) \to Ham(M,\omega) \tag 2.3
$$
the evaluation map $ev_1(\lambda) = \lambda(1) = \phi_H^1$.

\definition{Definition 2.2} Consider the metric
$d_H$ on $\PP(Ham(M,\omega),id)$ defined by
$$
d_H(\lambda,\mu) := \text{leng}(\lambda^{-1}\circ \mu), \quad
\lambda, \, \mu \in  \PP(Ham(M,\omega),id) \tag 2.4
$$
where $\lambda^{-1}\circ \mu$ is the Hamiltonian path $t \in [0,1]
\mapsto \lambda(t)^{-1}\mu(t)$. We call the induced topology  the
Hofer topology on $\PP(Ham(M,\omega),id)$. The Hofer topology on
$Ham(M,\omega)$ is the weakest topology for which the evaluation
map (2.3) is continuous.
\enddefinition
It is easy to see that this definition of the Hofer topology of
$Ham(M,\omega)$ coincides with the usual one induced by (1.1)
which also shows that the Hofer topology is meterizable. Of
course, nontriviality of the Hofer topology on $Ham(M,\omega)$ is
a difficult theorem which was proven by Hofer himself for $\C^n$
and by Lalonde and McDuff in its complete generality [LM1].

We say that two Hamiltonians $H$ and $K$ are {\it equivalent} if
their corresponding Hamiltonian paths $\phi_H$ and $\phi_K$ are
path-homotopic relative to the end points on $Ham(M,\omega)$. In
other words, for a given $\phi \in Ham(M,\omega)$ and $H$ and $K$
with
$$
H, \, K \mapsto \phi
$$
are equivalent if they are connected by one parameter family of
Hamiltonians $\{F^s\}_{0\leq s\leq 1}$ such that $F^0 = H, \, F^1
= K$ and
$$
F^s \mapsto \phi \quad \text{for all $s \in [0,1]$.}
$$
To emphasize the time-one map of the Hamiltonian path $\phi_H$, we
sometimes denote the Hamiltonian path $\phi_H$ also by $(\phi,H)$
and  by $h = [\phi, H] = [H]$ the equivalence class of $(\phi,H)$.
We denote by $\widetilde{Ham}(M,\omega)$ the set of equivalence
classes $[\phi,H]$ with the quotient topology on it. Although
$Ham(M,\omega)$ is not known to be locally path-connected in
general, $\widetilde{Ham}(M,\omega)$ can be considered as the
`universal covering space' in the \'etale sense (see [Oh7] for a
precise definition of topological \'etale covering).

One can easily check that $\overline H\# H'$ is given by the
formula
$$
(\overline H \# H')(t,x) = (H' - H)(t, (\phi_H^t)(x))
$$
and generates the flow $\phi_H^{-1}\circ \phi_{H'}$. Therefore we
can also write
$$
d_H(\phi_H,\phi_{H'}) = \text{leng}(\phi_H^{-1}\phi_{H'}) =
\|\overline H \# H' \|. \tag 2.5
$$

The following theorem from [Oh7] is an improvement of [Lemma 5.1,
Oh3]. For reader's convenience, we reproduce the proof from
[Oh3,7] here.

\proclaim{Theorem 2.3 [Theorem 5.1, Oh7]} The length minimizing
property of Hamiltonian path is closed in $\PP(Ham(M,\omega))$
with respect to the Hofer topology.
\endproclaim

\demo{Proof} ~ Suppose that \roster

\item the path $\phi_{G_i}$ converges to $\phi_{G_0}$
\item all $\phi_{G_i}$ are length minimizing in its homotopy class
relative to end points.
\endroster
Under these conditions, we need to prove that $\phi_{G_0}$ is
length minimizing in its homotopy class relative to the end points
again .

By definition of the Hofer topology, $\{G_i\}$ will satisfy $\|G_i
\#\overline G_0\| \to 0$ as $i \to \infty$. Suppose the contrary
that there exists $F$ such that $F\sim G_0$, but $\|F\| <
\|G_0\|$. Then there exists for some $\delta > 0$ such that
$$
\|F\|< \|G_0\| - \delta \tag 2.6
$$
Therefore we have
$$
\|F\| < \|G_i\| - {\delta \over 2} \tag 2.7
$$
for all sufficiently large $i$. We consider the Hamiltonian $F_i$
defined by
$$
F_i : = (G_i \# \overline{G_0})\# F \quad\text{i.e.,}
$$
This generates the flow $\phi^t_{G_i}\circ (\phi^t_{G_0})^{-1}
\circ \phi^t_F$ and so $F_i \sim G_i$. This implies, by the
hypothesis that $G_i$ are length minimizing over $[0,1]$, we have
$$
\|G_i\| \leq \|F_i\|. \tag 2.8
$$
On the other hand, we have
$$
\lim_{i \to \infty} \|F_i\|  = \lim_{i \to \infty} \|(G_i \#
\overline G_0) \# F)\| \leq \lim_{i \to \infty} (\|G_i \#
\overline G_0 \| + \| F\|) = \|F\| \tag 2.9
$$
Combining (2.7), (2.8) and (2.9), we get a contradiction. This
finishes the proof. \qed
\enddemo

We would like to note that in terms of the generating Hamiltonians
the Hofer topology is essentially $L^{(1,\infty)}$-topology which
is much weaker than e.g., $C^1$-topology. This point will play an
essential role when we will formulate our $C^1$ Perturbation
Conjecture on the Hamiltonian functions and reduce the Minimality
Conjecture to the $C^1$ Perturbation Conjecture in the very
strongly semi-positive case based on the main theorem of this
paper.

\head \bf \S 3. Chain level Floer theory and spectral invariants
\endhead

We first recall the construction of the spectral invariants from
[Oh5] briefly. Let $\Omega_0(M)$ be the set of contractible loops
and $\widetilde\Omega_0(M)$ be its standard covering space in the
Floer theory.  Note that the universal covering space of
$\Omega_0(M)$ can be described as the set of equivalence classes
of the pair $(\gamma, w)$ where $\gamma \in \Omega_0(M)$ and $w$
is a map from the unit disc $D=D^2$ to $M$ such that $w|_{\partial
D} = \gamma$: the equivalence relation to be used is that
$[\overline w \# w^\prime]$ is zero in $\pi_2(M)$. We say that
$(\gamma,w)$ is {\it $\Gamma$-equivalent} to $(\gamma,w^\prime)$
if and only if
$$
\omega([w'\# \overline w]) = 0 \quad \text{and }\, c_1([w\#
\overline w]) = 0 \tag 3.1
$$
where $\overline w$ is the map with opposite orientation on the
domain and $w'\# \overline w$ is the obvious glued sphere. And
$c_1$ denotes the first Chern class of $(M,\omega)$. We denote by
$[\gamma,w]$ the $\Gamma$-equivalence class of $(\gamma,w)$, by
$\widetilde\Omega_0(M)$ the set of $\Gamma$-equivalence classes
and by $\pi: \widetilde \Omega_0(M) \to \Omega_0(M)$ the canonical
projection. We also call $\widetilde \Omega_0(M)$ the
$\Gamma$-covering space of $\Omega_0(M)$. The unperturbed action
functional $\AA_0: \widetilde \Omega_0(M) \to \R$ is defined by
$$
\AA_0([\gamma,w]) = -\int w^*\omega. \tag 3.2
$$
Two $\Gamma$-equivalent pairs $(\gamma,w)$ and $(\gamma,w^\prime)$
have the same action and so the action is well-defined on
$\widetilde\Omega_0(M)$. When a periodic Hamiltonian $H:M \times
(\R/\Z) \to \R$ is given, we consider the functional $\AA_H:
\widetilde \Omega(M) \to \R$ defined by
$$
\AA_H([\gamma,w])= -\int w^*\omega - \int H(\gamma(t),t)dt
$$
We would like to note that {\it under this convention the maximum
and minimum are reversed when we compare the action functional
$\AA_G$ and the (quasi-autonomous) Hamiltonian $G$}.

We denote by $\text{Per}(H)$ the set of periodic orbits of $X_H$.
\medskip

\definition{Definition 3.1}  We define the {\it action spectrum}
of $H$, denoted as $\hbox{\rm Spec}(H) \subset \R$, by
$$
\hbox{\rm Spec}(H) := \{\AA_H(z,w)\in \R ~|~ [z,w] \in
\widetilde\Omega_0(M), z\in \text {Per}(H) \},
$$
i.e., the set of critical values of $\AA_H: \widetilde\Omega(M)
\to \R$. For each given $z \in \text {Per}(H)$, we denote
$$
\hbox{\rm Spec}(H;z) = \{\AA_H(z,w)\in \R ~|~ (z,w) \in
\pi^{-1}(z) \}.
$$
\enddefinition

Note that $\text {Spec}(H;z)$ is a principal homogeneous space
modelled by the period group of $(M,\omega)$
$$
\Gamma_\omega = \Gamma(M,\omega) := \{ \omega(A)~|~ A \in \pi_2(M)
\}
$$
and
$$
\hbox{\rm Spec}(H)= \cup_{z \in \text {Per}(H)}\text {Spec}(H;z).
$$
Recall that $\Gamma_\omega$ is either a discrete or a countable
dense subset of $\R$. It is trivial, i.e., $\Gamma_\omega = \{0\}$
in the weakly exact case. The followings were proved in [Oh3,4].

\proclaim{Lemma 3.2}~ $\hbox{\rm Spec}(H)$ is a measure zero
subset of $\R$.
\endproclaim

For given $\phi \in {\cal  H }am(M,\omega)$, we denote
$$
\HH_m(\phi) = \{ H ~|~ H \mapsto \phi, \, H\,  \text{normalized}
\}.
$$

\proclaim{Lemma 3.3} Let $F,\, G \in \HH_m(\phi)$ and $F\sim G$.
Then we have
$$
\text{\rm Spec}(G) = \text{\rm Spec}(F)
$$
as a subset of $\R$.
\endproclaim

This enables us to define the action spectrum over the universal
covering space $\widetilde{Ham}(M,\omega)$.

\definition{Definition 3.4} Let $h \in \widetilde{Ham}(M,\omega)$
and let $h = [\phi,H]$ for some Hamiltonian $H$ with
$\phi = \phi_H^1$. Then we define the action spectrum of $h$ by
$$
\text{\rm Spec}(h) = \text{\rm Spec}(H)
$$
for a (and so any) representative $[\phi,H]$.
\enddefinition

Next we  briefly recall the basic chain level operators in the
Floer theory, and the definition and basic properties of spectral
invariants $\rho(H;a)$ from [Oh5].

For each given generic time-periodic $H: M \times S^1 \to \R $, we
consider the free $\Q$ vector space over
$$
\text{Crit}\AA_H = \{[z,w]\in \widetilde\Omega_0(M) ~|~ z \in
\text{Per}(H)\}.
$$
To be able to define the Floer boundary operator correctly, we
need to complete this vector space downward with respect to the
real filtration provided by the action $\AA_H([z,w])$ of the
critical point $[z, w]$. More precisely, following [HS], [Oh3], we
introduce
\medskip

\definition{Definition 3.5} (1) We call the formal sum
$$
\beta = \sum _{[z, w] \in \text{Crit}\AA_H} a_{[z, w]} [z,w], \,
a_{[z,w]} \in \Q \tag 3.3
$$
a {\it Floer Novikov chain} if there are only finitely many
non-zero terms in the expression (3.2) above any given level of
the action. We denote by $CF(H)$ the set of Novikov chains. We
often simply call them {\it Floer chains}, especially when we do
not need to work on the covering space $\widetilde\Omega_0(M)$ as
in the weakly exact case.

(2) Two Floer chains $\alpha$ and $\alpha'$ are said to be {\it
homologous} to each other if they satisfy
$$
\alpha' = \alpha + \part_H(\gamma)
$$
for some Floer chain $\gamma$. We call $\beta$ a {\it Floer cycle}
if $\part \beta = 0$.

(3) Let $\beta$ be a Floer chain in $CF(H)$. We define and denote
the {\it level} of the chain $\beta$ by
$$
\lambda_H(\beta) =\max_{[z,w]} \{\AA_H([z,w]) ~|~a_{[z,w]}  \neq
0\, \text{in }\, (3.3) \}
$$
if $\beta \neq 0$, and just put $\lambda_H(0) = +\infty$ as usual.

(4) We say that $[z,w]$ is a {\it generator} of or {\it
contributes} to $\beta$ and denote
$$
[z,w] \in \beta
$$
if $a_{[z,w]} \neq 0$.
\enddefinition

Let $J = \{J_t\}_{0\leq t \leq 1}$ be a periodic family of
compatible almost complex structures on $(M,\omega)$.

For each given such periodic pair $(J, H)$, we define the boundary
operator
$$
\part: CF(H) \to CF(H)
$$
considering the perturbed Cauchy-Riemann equation
$$
\cases
\frac{\part u}{\part \tau} + J\Big(\frac{\part u}{\part t}
- X_H(u)\Big) = 0\\
\lim_{\tau \to -\infty}u(\tau) = z^-,  \lim_{\tau \to
\infty}u(\tau) = z^+ \\
\endcases
\tag 3.4
$$
This equation, when lifted to $\widetilde \Omega_0(M)$, defines
nothing but the {\it negative} gradient flow of $\AA_H$ with
respect to the $L^2$-metric on $\widetilde \Omega_0(M)$ induced by
the metrics $g_{J_t}: = \omega(\cdot, J_t\cdot)$ . For each given
$[z^-,w^-]$ and $[z^+,w^+]$, we define the moduli space
$$
\MM_{(J,H)}([z^-,w^-],[z^+,w^+])
$$
of solutions $u$ of (3.4) satisfying
$$
w^-\# u \sim w^+. \tag 3.5
$$
$\part$ has degree $-1$ and satisfies $\part\circ \part = 0$.

When we are given a family $(j, \HH)$ with $\HH = \{H^s\}_{0\leq s
\leq 1}$ and $j = \{J^s\}_{0\leq s \leq 1}$, the chain
homomorphism
$$
h_{(j,\HH)}: CF(H^0) \to CF(H^1)
$$
is defined by the non-autonomous equation
$$
\cases \frac{\part u}{\part \tau} +
J^{\rho_1(\tau)}\Big(\frac{\part u}{\part t}
- X_{H^{\rho_2(\tau)}}(u)\Big) = 0\\
\lim_{\tau \to -\infty}u(\tau) = z^-,  \lim_{\tau \to
\infty}u(\tau) = z^+
\endcases
\tag 3.6
$$
where $\rho_i, \, i= 1,2$ is functions of the type $\rho :\R \to
[0,1]$,
$$
\align
\rho(\tau) & = \cases 0 \, \quad \text {for $\tau \leq -R$}\\
                    1 \, \quad \text {for $\tau \geq R$}
                    \endcases \\
\rho^\prime(\tau) & \geq 0
\endalign
$$
for some $R > 0$. We denote by
$$
\MM^{(j,\HH)}([z^-,w^-],[z^+,w^+])
$$
or sometimes with $j$ suppressed the set of solutions of (3.6)
that satisfy (3.5). The chain map $h_{(j,\HH)}$ is defined
similarly as $\part$ using this moduli space instead.
$h_{(j,\HH)}$ has degree 0 and satisfies
$$
\part_{(J^1,H^1)} \circ h_{(j,\HH)} = h_{(j,\HH)} \circ
\part_{(J^0,H^0)}.
$$

The following identity can be proven by a straightforward
calculation, but has played a fundamental role in [Oh2-7] and will
also play an essential role in the proof of Theorem $\text{A}'$
and B.

\proclaim{Lemma 3.6} Let $H, K$ be any Hamiltonian not necessarily
non-degenerate and $j = \{J^s\}_{s \in [0,1]}$ be any given
homotopy and $\HH^{lin} = \{H^s\}_{0\leq s\leq 1}$ be the linear
homotopy $H^s = (1-s)H + sK$. Suppose that (3.5) has a solution
satisfying (3.6). Then we have the identity
$$
\align \AA_K([z^+,w^+]) & - \AA_H([z^-,w^-]) \\
& = - \int \Big|\dudtau \Big|_{J^{\rho_1(\tau)}}^2 -
\int_{-\infty}^\infty \rho_2'(\tau)\Big(K(t,u(\tau,t)) -
H(t,u(\tau,t))\Big) \, dt\,d\tau  \tag 3.7
\endalign
$$
\endproclaim

Now we recall the definition and some basic properties of spectral
invariant $\rho(H;a)$ from [Oh5]. We refer readers to [Oh5] for
the complete discussion on general properties of $\rho(H;a)$.

\proclaim{Definition \& Theorem 3.7 [Oh5]} Let $a \neq 0$ be a
given quantum cohomology class in $QH^*(M)$, and denote by
$a^\flat \in FH_*$ the  Floer homology class dual to $a$ in the
sense of [Oh5]. For any given Hamiltonian path $\lambda = \phi_H
\in \PP(Ham(M,\omega), id)$, we define
$$
\rho(\lambda;a): = \rho(H;a) = \inf_{\alpha \in \ker \part_H} = \{
\lambda_H(\alpha) \mid [\alpha] = a^\flat \}
$$
where $a^\flat$ is the dual to the quantum cohomology class $a$ in
the sense of [Oh5]. We call any of these {\it spectral
invariants}. The map $\rho_a: \lambda = \phi_H \mapsto \rho(H;a)$
defines a continuous function
$$
\rho_a=\rho(\cdot; a): C^\infty([0,1]\times M,\R) \to \R
$$
with respect to Hofer topology, and for two smooth functions $H
\sim K$ it satisfies
$$
\rho(H;a) = \rho(K;a)
$$
for all $a \in QH^*(M)$. In particular, for each given $h \in
\widetilde{\HH am}(M,\omega)$, the following definition is
well-defined:
$$
\rho(h;a) = \rho(H;a)
$$
any representative $[\phi,H] = h$.
\endproclaim

Now we focus on the invariant $\rho(h; 1)$ for $1 \in QH^*(M)$. We
first recall the following quantities
$$
\align E^-(h) & = \inf_{[\phi,H] = h}
E^-(H) \tag 3.9\\
E^+(h) & = \inf_{[\phi,H] = h} E^+(H) \tag 3.10
\endalign
$$
defined for smooth $h$. The following is an immediate consequence
of the proofs of similar inequalities from [Theorem II, Oh5].

\proclaim{Proposition 3.9} Let $(M,\omega)$ be arbitrary,
especially non-exact, closed symplectic manifold. For any $h \in
\widetilde{\HH am} (M,\omega)$, we have
$$
\rho(h;1) \leq E^-(h) \quad \rho(h^{-1};1) \leq E^+(h). \tag 3.11
$$
\endproclaim

The following theorem is an immediate consequence of Theorem 3.7
and Proposition 3.9 applied to the smooth case.

\proclaim{Theorem 3.10} Let $G: [0,1] \times M \to \R$ be a
quasi-autonomous Hamiltonian. Suppose that $G$ satisfies
$$
\rho(G;1) = E^-(G) \tag 3.12
$$
Then $G$ is negative Hofer-length minimizing in its homotopy class
with fixed ends.
\endproclaim

So far in this section, we have presumed that the Hamiltonians are
time one-periodic. Now we explain how to dispose the periodicity
and extend the definition of $\rho(H;a)$ for arbitrary time
dependent Hamiltonians $H: [0,1] \times M \to \R$. Note that it is
obvious that the semi-norms $E^\pm(H)$ and $\|H\|$ are defined
without assuming the periodicity. For this purpose, the following
lemma from [Oh3] is important. We leave its proof to readers or to
[Oh3].

\proclaim{Lemma 3.11} Let $H$ be a given Hamiltonian $H : [0,1]
\times M \to \R$ and $\phi = \phi_H^1$ be its time-one map. Then
we can re-parameterize $\phi_H^t$ in time so that the
re-parameterized Hamiltonian $H'$ satisfies the following
properties: \roster
\item $\phi_{H'}^1 = \phi_H^1$
\item $H' \equiv 0$ near $t = 0, \, 1$ and in particular $H'$ is
time periodic. We call such Hamiltonians {\it boundary flat}.
\item Both $E^\pm(H' - H)$ can be made as small as we want
\item If $H$ is quasi-autonomous, then so is $H'$
\item For the Hamiltonians $H', \, H''$ generating any two such
re-parameterizations of $\phi_H^t$, there is canonical one-one
correspondences between $\text{Per}(H')$ and $\text{Per}(H'')$,
and $\text{Crit }\AA_{H'}$ and $\text{Crit }\AA_{H''}$ with their
actions fixed .
\endroster
Furthermore this re-parameterization is canonical with the
``smallness'' in (3) can be chosen uniformly over $H$ depending
only on the $C^0$-norm of $H$. In particular, this approximation
can be done with respect to the Hofer topology.
\endproclaim
In fact, the above approximation can be done under the stronger
topology, the {\it strong Hamiltonian topology} introduced by the
author [Oh7]. Since this stronger statement will not be needed in
the present paper, we will be content with stating this lemma
under the Hofer topology here.

\head \bf \S 4. Canonical fundamental Floer cycles
\endhead

In this section, we start with our study of length minimizing
property of Hamiltonian paths following the scheme used in [Oh6].
This section is largely a duplication of section 4 of [Oh6].
Partly because the details of the proofs, not the theorems
therein, will be needed in the proof of Theorem B in the next
section and also to make this paper self-contained, we repeat the
whole section 4 of [Oh6] here with minor simplification for the
readers' convenience.

We first recall the basic definitions in relation to the dynamics
of Hamiltonian flows.

\definition{\bf Definition 4.1} Let $H: M \times [0,1] \to \R$ be a
Hamiltonian which is not necessarily time-periodic and $\phi_H^t$
be its Hamiltonian flow. \par

\roster \item We call a point $p\in M$ a {\it time $T$ periodic
point} if $\phi_H^T(p)=p$. We call $t \in [0,T] \mapsto
\phi_H^t(p)$ a {\it contractible time $T$-periodic orbit} if it is
contractible. \par

\item When $H$ has a fixed critical point $p$ over $t \in
[0,T]$, we call $p$ {\it over-twisted} as a time $T$-periodic
orbit if its linearized flow $d\phi_H^t(p); \, t\in [0,T]$ on
$T_pM$ has a closed trajectory, other than the fixed origin, of
period less than or equal to $T$. Otherwise we call it {\it
under-twisted}. If in addition the linearized flow has
only the origin as the fixed point, then we call {\it generically
under-twisted}.
\endroster
\enddefinition

For the proof of Theorem B, we need to unravel the definition of
the particular spectral invariant $\rho(G;1)$ from [Oh5], and its
realization as the level of some optimal Floer cycles for
nondegenerate (one periodic) Hamiltonians $G$. According to the
definition from [Oh5] for nondegenerate Hamiltonians, we consider
the Floer homology class dual to the quantum cohomology class $1
\in H^*(M) \subset QH^*(M)$, which we denote by $1^\flat$
following the notation of [Oh5] and call the {\it semi-infinite}
fundamental class of $M$. Then according to [Definition 4.2 \&
Theorem 4.5, Oh5], we have
$$
\rho(G;1) = \inf_{\gamma}\{\lambda_G(\gamma) \mid \gamma \in
\ker\part_G \subset CF(G)\, \text{with }\, [\gamma] = 1^\flat \}.
\tag 4.1
$$
$\rho$ is then extended to arbitrary Hamiltonians by continuity in
$C^0$-topology. Therefore to prove Theorem B for the nondegenerate
Hamiltonians, we need to first construct cycles $\gamma$ with
$[\gamma] = 1^\flat$ whose level $\lambda_G(\gamma)$ become
arbitrarily close to $E^-(G)$, and then to prove that the cycle
cannot be pushed down by the Cauchy-Riemann flow.

We recall the following concept of homological essentialness in
the chain level theory from [Oh3,6].

\definition{Definition 4.2} We call a Floer cycle $\alpha \in CF(H)$
{\it tight} if it satisfies the following non-pushing down
property under the Cauchy-Riemann flow (3.4): for any Floer cycle
$\alpha' \in CF(H)$ homologous to $\alpha$ (in the sense of
Definition 3.1 (2)), it satisfies
$$
\lambda_H(\alpha') \geq \lambda_H(\alpha). \tag 4.2
$$
\enddefinition

Now we will need to construct a {\it tight} fundamental Floer
cycle of nondegenerate quasi-autonomous $H$ whose level is
precisely $E^-(H)$. As a first step, we first construct a
fundamental cycle of $H$ whose level is $E^-(H)$ which may not be
tight in general.

We choose a Morse function $f$ such that $f$ has the unique global
minimum point $x^-$ and
$$
f(x^-) = 0, \quad f(x^-) < f(x_j) \tag 4.3
$$
for all other critical points $x_j$. Then we choose a fundamental
Morse cycle
$$
\alpha =\alpha_{\e f} = [x^-,w_{x^-}] + \sum_j a_j [x_j,w_{x_j}]
$$
as in [Oh3] where $x_j \in \text{Crit }_{2n} (-f)$. Recall that
the {\it positive} Morse gradient flow of $\e f$ corresponds to
the {\it negative} gradient flow of $\AA_{\e f}$ in our
convention.

Considering Floer's homotopy map $h_\LL$ over the linear path
$$
\LL: \, s\mapsto (1-s) \e f +s H
$$
for sufficiently small $\e > 0$, we transfer the above fundamental
Morse cycle $\alpha$ and define a fundamental Floer cycle of $H$
by
$$
\alpha_H: = h_{\LL}(\alpha) \in CF(H). \tag 4.4
$$
We call this particular cycle {\it the canonical fundamental Floer
cycle} of $H$.

The following important property of this fundamental cycle was
proved by Kerman-Lalonde [KL] for the aspherical case and by the
author [Oh6] in general.

\proclaim{Proposition 4.3} Suppose that  $H$ is a generic
one-periodic Hamiltonian such that $H_t$ has the unique
non-degenerate global minimum $x^-$ which is fixed and generically
under-twisted for any $t \in [0,1]$. Suppose that $f: M \to \R$ is
a Morse function such that $f$ has the unique global minimum point
$x^-$ and $f(x^-)=0$. Then the canonical fundamental cycle has the
expression
$$
\alpha_H = [x^-, w_{x^-}] + \beta \in CF(H)\tag 4.5
$$
for some Floer Novikov chain $\beta \in CF(H)$ with the inequality
$$
\lambda_H(\beta) < \lambda_H([x^-,w_{x^-}]) = \int_0^1 - H(t,
x^-)\, dt. \tag 4.6
$$
In particular its level satisfies
$$
\align \lambda_H(\alpha_H) & = \lambda_H([x^-,w_{x^-}])\tag 4.7\\
& = \int_0^1 - H(t, x^-)\, dt = \int_0^1 -\min H\, dt.
\endalign
$$
\endproclaim

The proof is based on the following simple fact used by Kerman and
Lalonde (see the proof of [Proposition 4.2, KL]).  We refer to
[Lemma 3.5, Oh6] for the details of its proof.

\proclaim{Lemma 4.4} Let $H$ and $f$ as in Proposition 3.3. Then
for all sufficiently small $\e > 0$, the function $G^H$ defined by
$$
G^H(t,x) = H(t, x^-) + \e f
$$
satisfies
$$
\aligned G^H(t,x^-) & = H(t, x^-) \\
G^H(t,x)& \leq H(t,x)
\endaligned
\tag 4.8
$$
for all $(t,x)$ and equality holds only at $x^-$.
\endproclaim

\demo{Proof of Proposition 4.3} Since $x^-$ is a under-twisted
fixed minimum of both $H$ and $f$, we have the Conley-Zehnder
index
$$
\mu_H([x^-, w_{x^-}])  = \mu_{\e f}([x^-,w_{x^-}]) ( = n)
$$
and so the moduli space $\MM^\LL([x^-,w_{x^-}],[x^-,w_{x^-}])$ has
dimension zero. Let $u \in \MM^\LL([x^-,w_{x^-}],[x^-,w_{x^-}])$.

We note that the Floer continuity equation (3.6) for the linear
homotopy
$$
\LL: s \to (1-s) \e f + s H
$$
is unchanged even if we replace the homotopy by the homotopy
$$
\LL': s \to (1-s) G^H + s H.
$$
This is because the added term $H(t,x^-)$ in $G^H$ to $\e f$ does
not depend on $x \in M$ and so
$$
X_{\e f} \equiv X_{G^H}.
$$
Therefore $u$ is also a solution for the continuity equation (3.6)
under the linear homotopy $\LL'$. Using this,  we derive the
identity
$$
\aligned \int \Big|\dudtau \Big|_{J^{\rho_1(\tau)}}^2\, dt\,d\tau
& =
\AA_{G^H}([x^-,w_{x^-}]) - \AA_H([x^-,w_{x^-}]) \\
\quad & - \int_{-\infty}^\infty \rho'(\tau)\Big(H(t,u(\tau,t))\,
dt\, d\tau - G^H(t,u(\tau,t))\Big) \, dt\,d\tau
\endaligned
\tag 4.9
$$
from (3.7). Since we have
$$
\AA_H([x^-,w_{x^-}]) =\AA_{G^H}([x^-,w_{x^-}]) = \int_0^1 - \min
H\, dt \tag 4.10
$$
and $G^H \leq H$, the right hand side of (4.9) is non-positive.
Therefore we derive that $\MM^\LL([x^-,w_{x^-}],[x^-,w_{x^-}])$
consists only of the constant solution $u \equiv x^-$. This in
particular gives rise to the matrix coefficient of $h_\LL$
satisfying
$$
\langle [x^-,w_{x^-}], h_{\LL}([x^-,w_{x^-}])\rangle = \#
(\MM^\LL([x^-,w_{x^-}],[x^-,w_{x^-}])) = 1.
$$
Now consider any other generator of $\alpha_H$
$$
[z,w] \in \alpha_H \quad \text{with }\,  [z,w] \neq [x^-,w_{x^-}].
$$
By the definitions of $h_\LL$ and $\alpha_H$, there is a generator
$[x,w_x] \in \alpha$ such that
$$
\MM^\LL([x,w_x],[z,w]) \neq \emptyset. \tag 4.11
$$
Then for any $u \in \MM^\LL([x,w_x],[z,w])$, we have the identity
from (3.7)
$$
\aligned  \AA_H([z,w]) - & \AA_{G^H}([x,w_x]) = -\int
\Big|\dudtau \Big|_{J^{\rho_1(\tau)}}^2\, dt\, d\tau \\
& \quad - \int_{-\infty}^\infty \rho'(\tau)\Big(H(t,u(\tau,t)) -
G^H(t,u(\tau,t))\Big) \, dt\,d\tau. \endaligned
$$
Since $-\int \Big|\dudtau \Big|_{J^{\rho_1(\tau)}}^2 \leq 0$, and
$G^H \leq H$,  we have
$$
\AA_H([z,w]) \leq \AA_{G^H}([x,w_x])\tag 4.12
$$
with equality holding only when $u$ is stationary. There are two
cases to consider, one for the case of $x = x^-$ and the other for
$x = x_j$ for $x_j \neq x^-$ for $[x_j,w_{x_j}] \in \alpha$.

For the first case, {\it since we assume $[z,w] \neq
[x^-,w_{x^-}]$}, $u$ cannot be constant and so the strict
inequality holds in (4.12), i.e,
$$
\AA_H([z,w]) < \AA_{G^H}([x^-,w_{x^-}]). \tag 4.13
$$
For the second case, we have the inequality
$$
\AA_H([z,w]) \leq \AA_{G^H}([x_j,w_{x_j}])\tag 4.14
$$
for some $x_j \neq x^-$ with $[x_j,w_{x_j}] \in \alpha$. We note
that (4.3) is equivalent to
$$
\AA_{G^H}([x_j,w_{x_j}]) < \AA_{G^H}([x^-,w_{x^-}]).
$$
This together with (4.14) again give rise to (4.13). On the other
hand we also have
$$
\AA_{G^H}([x^-,w_{x^-}]) = \AA_H([x^-,w_{x^-}])
$$
because $G^H(t,x^-) = H(t,x^-)$ from (4.5). Altogether, we have
proved
$$
\AA_H([z,w]) < \AA_H([x^-,w_{x^-}]) = \int_0^1 -H(t,x^-) \, dt
$$
for any $[z,w] \in \alpha_H$ with $[z,w] \neq [x^-,w_{x^-}]$. This
finishes the proof of (4.3). \qed\enddemo

\head \bf\S 5. Proof of Theorem $\text{A}'$ and Theorem B
\endhead

In this section, we will give the proof of Theorem B and hence
Theorem $\text{A}'$ via the scheme of the previous section. We now
rephrase Theorem B here.

\proclaim{Theorem 5.1} Suppose that $(M,\omega)$ is very strongly
semi-positive. Let $H, \, H'$ be as in Theorem $\text{A}'$. Then
we have
$$
\rho(H';1) = E^-(H') = \int_0^1 -\min H' \, dt. \tag 5.1
$$
\endproclaim

The theorem is an immediate consequence of the following tightness
result of the canonical fundamental Floer cycle constructed in
section 4. Here the very strongly semi-positive condition enters
in an essential way.

\proclaim{Proposition 5.2} Let $(M,\omega)$ and $H, \, H'$ be as
in Theorem 5.1 and let $\alpha_{H'}$ be the canonical fundamental
Floer cycle of $H'$. Then $\alpha_{H'}$ is tight: i.e., for any
Floer Novikov cycle $\alpha \in CF(H')$ homologous to
$\alpha_{H'}$, we have
$$
\lambda_{H'}(\alpha) \geq \lambda_{H'}(\alpha_{H'}). \tag 5.2
$$
In particular, we have
$$
\rho(H';1) = \lambda_{H'}(\alpha_{H'}) \, (= E^-(H')).
$$
\endproclaim
\demo{Proof} Suppose that $\alpha$ is homologous to $\alpha_{H'}$,
i.e.,
$$
\alpha = \alpha_{H'} + \part_{H'}(\gamma) \tag 5.3
$$
for some Floer Novikov chain $\gamma \in CF(H')$. We need to check
whether $[x^-,w_{x^-}]$ in the cycle $\alpha_{H'}$ can be
cancelled by $\part_{H'}(\gamma)$ for a suitable choice of a Floer
chain $\gamma$. For this purpose, we study the matrix element
$$
\langle [x^-, w_{x^-}], \part_{H'} [z,w] \rangle
$$
for each $[z,w] \in \text{Crit}\AA_{H'}$, which in turn need to
study the integers
$$
\#(\MM_{J,H'}([z,w],[x^-,w_{x^-}])).
$$
There are two types of $[z,w]$ to consider, one $z = x^-$ and the
other $z \neq x^-$.

We first state the following perturbation lemma whose proof we
omit and refer to [Appendix, KL] for the details, which however
considers only the $C^2$-perturbations, but immediately
generalizes to the $C^1$-perturbations {\it as long as we do not
perturb $H$ near $x^-$}. See the proof of Proposition 6.4 for some
relevant adjustment needed.

\proclaim{Lemma 5.3} Suppose that $H$ has a generically
under-twisted local minimum point $x^-$. Then there exists a fixed
neighborhood $U_{x^-}$ depending only on $H$ such that \roster
\item for any sufficiently $C^1$-small
perturbation $H'$ of $H$ with $H' \equiv H$ on $U_{x^-}$, $x^-$ is
the only critical point of $H'$ in $U_{x^-}$ which is generically
under-twisted and
\item
for any non-constant contractible periodic orbit $z$ of
$H^\prime$, $\text{Im }z \subset M \setminus U_{x^-}$.
\endroster
\endproclaim

We also need one more simple lemma about the lower bound for the
energy of the Floer trajectory connecting $x^-$ and any other
periodic orbits $z$. Again we omit its proof which is a simple
consequence of compactness arguments (see e.g., [Oh3] for such an
argument).

\proclaim{Lemma 5.4} Let $H$ be any smooth Hamiltonian which {\it
is not necessarily regular}. Suppose $H$ has the unique critical
point $x$ in a neighborhood $U_x$ and the image of no periodic
orbits intersect $U_x$. Then there exists a constant $ e_{H}> 0$
such that for any finite energy solution $u: \R \times S^1 \to M$
of (3.4) such that
$$
u(-\infty) = z, \, u(\infty)=x
$$
with $z \neq x$, we have
$$
\int\Big|\frac{\part u}{\part \tau}\Big|^2_J \geq e_H > 0\tag 5.4
$$
where $e_H$ does not depend on $u$.
\endproclaim

Let $U_{x^-}$ be as in Lemma 5.3. We recall that only the critical
points $[z,w]$ with
$$
\mu_{H'}([z,w])= \mu_{H'}([x^-,w_{x^-}]) + 1 = n +1 \tag 5.5
$$
can have non-zero matrix element $\langle [x^-,w_{x^-}],
\partial_{H'} [z,w]\rangle$.
Here the latter equality follows from the
generically-undertwistedness of $[x^-,w_{x^-}]$ and and the
definition of the Conley-Zehnder index. To provide some intuition
on why this is so, we quote the following general result
concerning the Conley-Zehnder index and the eigenvalues of the
linearization matrix of the autonomous Hamiltonians. We would like
to note that [SZ] uses a different convention of the Hamiltonian
vector field from ours. According to their convention, our $X_H$
is the negative of theirs. For example, if we apply
$S=-d^2G(x^-)$, we have $\mu_G([x^-,w_{x^-}]) = 2n -n =n$ since
$x^-$ is a maximum point of $-G$. One may regard the condition of
generically undertwistedness in Definition 4.1 [KL] as the
non-autonomous analog to this case.

\proclaim{Lemma 5.5 [Theorem 3.3 (iv), SZ]} Consider the matrix
$\Psi(t) = \exp{(J_0St)}$ where $J_0$ is the standard almost
complex structure on $\R^{2n}\cong\C^n$ and $S = S^T \in
M^{2n\times 2n}(\R)$ is a non-singular symmetric matrix such all
the eigenvalues $\lambda$ satisfies
$$
|\lambda| < 2 \pi. \tag 5.6
$$
Then the Conley-Zehnder index $\mu(\Psi)$ satisfies
$$
\mu(\Psi) = \mu^-(S) - n \tag 5.7
$$
where $\mu^-(S)$ denotes the number of negative eigenvalues of $S$
counted with multiplicity.
\endproclaim

We also have the general index formula for the Conley-Zehnder
index and the first Chern number
$$
\mu_{H'}([z,w]) = \mu_{H'}([z,w']) + 2 c_1([w\#\overline w']).
\tag 5.8
$$
(In the literature there are various other different conventions
used in relation to the definition of $X_H$, that of action
functional, and also both homological and cohomological notations
have been used. This makes the question about the correct sign for
the corresponding formula very confusing. All of our conventions
are the same with those in [En] and [Oh3], {\it except that there
is an error in the corresponding formula: in [En] and [Oh3], this
formula is written as
$$
\mu_{H'}([z,w]) = \mu_{H'}([z,w']) - 2 c_1([w\#\overline w']).
$$
This is incorrect.})

At this stage, we first consider the case where $z = x^-$, i.e.,
$[z,w] = [x^-,w]$ for an arbitrary $w$ bounding disc of $x^-$
which is nothing but a sphere passing through the point $x^-$. In
this case (5.8) becomes
$$
\mu_{H'}([x^-,w]) = \mu_{H'}([x^-,w_{x^-}]) + 2 c_1([w\#\overline
w_{x^-}])
$$
when applied to $[z,w'] = [x^-,w_{x^-}]$. We also recall that,
only when
$$
\mu_{H'}([x^-,w]) - \mu_{H'}([x^-,w_{x^-}]) = 1,
$$
there can be a Floer trajectory that is issued at $[x^-,w]$ and
landing at $[x^-,w_{x^-}]$, and provides a non-trivial matrix
element. However this formula is impossible because $2c_1(w \#
\overline w_{x^-})$ is an even number. Therefore the matrix
elements $\langle [x^-,w_{x^-}],\partial_{H'} [z,w]\rangle$ are
all zero when $z = x^-$. {\it For this we do not need the
condition of $\e$-positively $\mu$-undertwistedness of $H'$ but
only the generically undertwistedness of $x^-$.}

Next we consider the case $z \neq x^-$. We first recall that
$\mu_{H'}([x^-,w_{x^-}]) = n $ and so for $[z,w]$ to give a
nontrivial contribution to the Floer matrix element $\langle
[x^-,w_{x^-}],\part_{H'} [z,w] \rangle$ it must be the case that
$$
\mu_{H'}([z,w])= n + 1. \tag 5.9
$$
By the hypothesis of the $\mu$-undertwistedness on $z$ in Theorem
5.1, we have a bounding disc $w_z$ such that the Conley-Zehnder
index of $[z,w_z]$ satisfies
$$
-n \leq \mu_{H'}([z,w_z]) \leq n \tag 5.10
$$
for any one-periodic orbit $z$ of $H'$.

Now suppose that there exists $[z,w] \in \text{Crit}\AA_{H'}$ that
gives a non-trivial matrix coefficient
$$
\langle [x^-,w_{x^-}],\partial_{H'} [z,w]\rangle \neq 0
$$
and let $u$ be a trajectory from $[z,w]$ to $[x^-,w_{x^-}]$. We
consider the glued sphere $w_z \# u\# \overline{w}_{x^-} = w_z \#
u$ and $w\# \overline{w}_{z}$.

Combining (5.8)-(5.10) applied to $w'=w_z$, we derive
$$
1 \leq 2c_1([w \# \overline w_z] \leq 2n+1
$$
and hence
$$
\frac{1}{2} \leq c_1([w \# \overline w_z]) \leq n+ \frac{1}{2}.
$$
Because $c_1$ is an integer, we indeed must have
$$
0 < c_1([w \# \overline w_z]) \leq n
$$
or equivalently
$$
-n \leq c_1([\overline w \# w_z]) < 0. \tag 5.11
$$
({\it This inequality is the origin of our imposing the very
strongly semi-positive condition in Theorem A.})

Next we consider the actions of $[z,w]$. From (3.5),  $w\#
u\#\overline{w}_{x^-} = w\#u$ is contractible, and so we have
$$
0 = \int_{w\# u}\omega = \int_{w\# \overline{w}_z} \omega +
\int_{w_z\# u}\omega \tag 5.12
$$
We rewrite the second term into
$$
\int_{w_z\# u}\omega = \int w_z^*\omega + \int u^*\omega. \tag
5.13
$$
For the second term here, we have
$$
\align \int u^*\omega & = \int_{-\infty}^\infty \int_0^1
\omega\Big(
\frac{\part u}{\part \tau}, \frac{\part u}{\part t}\Big) \, dt\, d\tau \\
& =\int_{-\infty}^\infty \int_0^1 \omega\Big( \frac{\part u}{\part
\tau}, J \frac{\part u}{\part \tau} + X_{H'}(u)\Big)\, dt\, d\tau \\
& = \int \Big|\frac{\part u}{\part \tau}\Big|^2_J +
\int_{-\infty}^\infty \int_0^1 \omega\Big(\frac{\part u}{\part
\tau}, X_{H'}(u)\Big)\, dt\, d\tau \\
& = \int \Big|\frac{\part u}{\part \tau}\Big|^2_J + \int_0^1
\Big(- H'(t,x^-) + H'(t,z(t)\Big)\, dt \geq \int \Big|\frac{\part
u}{\part \tau}\Big|^2_J. \tag 5.14
\endalign
$$
For the last inequality, we have used the fact that $x^-$ is the
(global) minimum point of $H'$.

On the other hand, since $z$ does not intersect $U_{x^-}$, it
follows from Lemma 5.4 that we have
$$
\int \Big| \frac{\part u}{\part \tau}\Big|^2_J \geq e_H >0. \tag
5.15
$$
Therefore, combining (5.12)-(5.15), we obtain
$$
\int_{w_z\# u}\omega = \int w_z^*\omega + \int u^*\omega \geq \int
w_z^*\omega + e_H.
$$
And  by the condition on $w_z$ in the definition of
$\e$-positively $\mu$-undertwistedness of $H'$, we must have
$$
\int w_z^*\omega > - \e, \quad 0 < \e < e_H
$$
and so
$$
\int_{w_z\# u}\omega > 0 \tag 5.16
$$
Then (5.12) and (5.16) imply that we have $\int_{w\#
\overline{w}_z} w < 0$, i.e.,
$$
\int_{\overline w\# w_z} w > 0 \tag 5.17
$$
for all critical points $[z,w]$ such that $z \neq x^-$ and
$$
\langle [x^-,w_{x^-}], \part_{H'} [z,w] \rangle \neq 0.
$$
However it follows from the very strongly semi-positive hypothesis
that (5.11) and (5.17) together are prohibited and so there is no
critical point $[z,w], \, z \neq x^-$ such that $\langle
[x^-,w_{x^-}],
\part_{H'} [z,w] \rangle \neq 0$.

Altogether, addition of $\part_{H}(\gamma)$ to $\alpha_{H'}$
cannot kill the term $[x^-,w_{x^-}]$  away from the cycle
$$
\alpha_{H'} = [x^-, w_{x^-}] + \beta
$$
in (5.3), and hence we have
$$
\lambda_{H'}(\alpha)= \lambda_{H'}(\alpha_{H'} +
\part_{H'}(\gamma)) \geq \lambda_{H'}([x^-,w_{x^-}]) \tag 5.18
$$
by the definition of the level $\lambda_{H'}$.  Combining (4.10)
and (5.18), we have finished the proof (5.2). \qed\enddemo

\definition{Remark 5.7}
Note that (5.17) is not possible for the weakly exact case and so
already proves Theorem B in that case {\it without assuming the
$\e$-positively $\mu$-undertwistedness hypotheses for $H'$}.
\enddefinition

\demo{Finish-up of the proof of Theorem $\text{A}'$} Similar
consideration simultaneously applied to $x^+$, then proves
$$
\rho(\overline{H'};1) = E^-(\overline{H'}).
$$
Now let $F \sim H'$. Then Theorem 3.7 implies 
$$
\rho(F;1) = \rho(H';1).
$$
Combining this with (3.11) and Theorem 5.1 for $H'$ and
$\overline H'$, we derive 
$$
\align E^-(H') &\leq E^-(F) \tag 5.19\\
E^-(\overline{H'}) & \leq E^-(\overline F). \tag 5.20
\endalign
$$
Recall that (5.20) is equivalent to
$$
E^+(H') \leq E^+(F). \tag 5.21
$$
Adding (5.19) and (5.21) then gives rise to the inequality $\|H'\|
\leq \|F\|$. This proves $H'$ is length minimizing. Now let $\e
\to 0$. Then we have $H' \to H$ in the $C^1$-topology and so
$\phi_{H'} \to \phi_H$ in the Hofer topology. Since
Proposition 5.2 proves that all $\phi_{H'}:t \mapsto \phi_{H'}^t$
is length minimizing, by letting $\e \to 0$, we derive from
Theorem 2.3 that the Hamiltonian path $\phi_H: t \mapsto \phi_H^t$
itself is length minimizing. \qed\enddemo

\head \bf \S6. Minimality and the $C^1$ Perturbation Conjecture
\endhead

In this section, we restrict to the case of autonomous
Hamiltonians $G$ and discuss how one might try to prove the
Minimality Conjecture, at least for the very strongly
semi-positive $(M,\omega)$.

The following perturbation conjecture seems to be the crux towards
the proof of the Minimality Conjecture at least for the very
strongly semi-positive cases. Applying Lemma 5.3 and Lemma 5.4 to
the autonomous Hamiltonian $G$ in the Minimality Conjecture, the
conjecture for the very strongly semi-positive case will
immediately follow from the following general perturbation
conjecture, via Theorem $\text{A}'$ and B.

\proclaim{[$C^1$ Perturbation Conjecture]} Let $(M,\omega)$ be any
symplectic manifold and let $G$ an autonomous Hamiltonian as in
the Minimality Conjecture and let $\e
>0$ be given. Then we can choose a quasi-autonomous smooth Hamiltonian
$H$ that satisfies the following properties in addition: \roster

\item $H$ is non-degenerate in the Floer
theoretic sense and $\e$-positively $\mu$-undertwisted.

\item $H$ has the unique fixed minimum point $x^-$ and the
unique fixed maximum point $x^+$ that are generically under-twisted

\item The $C^1$-norm $\|G - H\|_{C^1}$ can be made
arbitrarily small as we want.
\endroster
\endproclaim
A more conservative version of the conjecture will be the one
restricted to the case when the autonomous Hamiltonian has only
{\it isolated} critical points.

The only non-trivial point of this conjecture lies in the
condition `$\e$-positively $\mu$-undertwistedness'. Furthermore we
will now show that $\e$-positivity can be automatically ensured by
any sufficiently $C^1$-small perturbation of $G$. For the
remaining section, we will explain this claim and how the $C^1$
Perturbation Conjecture implies the Minimality Conjecture on the
very strongly semi-positive symplectic manifold.

We first introduce the notion of {\it center of mass} of the
closed curve

\proclaim{Definition \& Lemma 6.1} Fix any Riemannian metric e.g.,
the already used compatible metric $g =\omega(\cdot, J\cdot)$ on
$M$. We denote by $\text{\rm inj}(g)>0$ the injectivity radius and
by $\exp$ the exponential map of $g$. Let $z: S^1 \to M$ be a
closed continuous curve whose diameter is less than the
injectivity radius of $M$. Then there exists a unique point $x_z$
and a unique closed continuous curve
$$
\xi_z: S^1 \to T_{x_z}M
$$
such that
$$
z(t) = \exp_{x_z}\xi_z(t). \tag 6.1
$$
The point $x_z$ is called the {\it center of mass of $z$}. If $z$
is $C^k$, then so is $\xi_z$.
\endproclaim
We refer readers to [K], for example, for the proof of this lemma.

Next we introduce the notion of {\it canonical small bounding
disc} of any closed curve $z: S^1 \to M$ that is contained in a
{\it convex} ball in the sense of Riemannian geometry.

\definition{Definition 6.2} Let $B\subset M$ be a closed strongly convex
ball. Supposed that the image of a continuous closed curve $z: S^1
\to M$ is contained in $Int B$. Consider the center of mass
representation (6.1) of the curve $z$. We call the disc $w_z:
D=D^2 \to M$ defined by
$$
w_z(r,t) = \exp_{x_z}(r\xi_z(t)) \tag 6.2
$$
in the polar coordinates the {\it canonical small bounding disc}
of $z$. Furthermore this representation is independent of
reparameterization of the curve.
\enddefinition

It follows that $w_z$ is always continuous, and is unique up to
homotopy in that any two bounding discs of $z$ whose images are
contained in the convex ball $B$ are homotopic to each other
relative to the boundary.

It is also $C^k$ away from $0$ if $z$ is $C^k$. In particular if
$z$ is $C^1$, the integral
$$
\int_{w_z} \omega
$$
is well-defined by defining it to be the integral of any
$C^1$-approximation of $w_z$ inside $Int B$ and fixing its
boundary to be $z$. With this definition, we have the following
easy lemma.

\proclaim{Lemma 6.3} Consider the compatible metric $g =
\omega(\cdot, J\cdot)$. If $z$ is as in Definition 6.2 and is
$C^1$, then the Riemannian area $Area_g(w_z)$ is well-defined and
we have the inequality
$$
\Big|\int_{w_z} \omega \Big| \leq Area_g(w_z). \tag 6.3
$$
\endproclaim
With these preparations, we prove the following proposition.

\proclaim{Proposition 6.4} Let $G$ be the autonomous Hamiltonian
as in the Minimality Conjecture. Then for any given small $0 < \e
< \frac{1}{2}\text{\rm inj}(g)$, there exists sufficiently small
$\delta_1, \, \delta_2 >0$ such that for any smooth Hamiltonian
$H:[0,1] \times M \to \R$ with $\|H-G\|_{C^1}< \delta_1$ the
following holds: for any periodic orbit $z$ of period $T$ with
$1-\delta_2 \leq T \leq 1+\delta_2$, we have \roster
\item
$$
\text{\rm diam}(z) \leq \frac{1}{2}\text{\rm inj}(g) \tag 6.4
$$
and
$$
d_{C^1}(z,x_z) \leq \e \tag 6.5
$$
\item and the canonical small bounding disc $w_z:D \to M$ with $\part
w_z = z$ satisfies
$$
\int_D w_z^*\omega \geq - \e. \tag 6.6
$$
\endroster
\endproclaim
\demo{Proof} We will prove this by contradiction. Suppose the
contrary that there exists some $\e_0 >0$ for which we can choose
sequences of one-periodic Hamiltonians $H_i$ with $\|H_i -
G\|_{C^1} \to 0$ and a sequence $z_i$ of periodic orbit of
$X_{H_i}$ of period $|T_i - 1| \to 0$ such that one of the
following occurs
\medskip

(i) either $\text{diam}(z_i) > \frac{1}{2}\text{inj}(g)$ or
\smallskip

(ii) $\text{diam}(z_i) \leq \frac{1}{2}\text{inj}(g)$ but
$d_{C^1}(z_i,x_{z_i}) > \e_0$ or
\smallskip

(iii) $\text{diam}(z_i) \leq \frac{1}{2}\text{inj}(g)$ and
$d_{C^1}(z_i,x_{z_i}) \leq \e_0$ but $\int_D w_z^*\omega < -
\e_0$.
\medskip

Using the fact that the period $T_i$ satisfies $|T_i - 1| \to 0$,
we consider the repamameterized  Hamiltonians $H_i'$ defined by
$$
H_i'(t,x): = T_i H_i(T_i t,x).
$$
We note that $H_i'$ is {\it not} one-periodic but $(1/T_i)$-periodic.
But this is irrelevant  as long as it is a smooth function on $\R
\times M$ which certainly is. The repamameterized  orbit
$$
\widetilde z_i(t) = z_i(T_it)
$$
now defines a smooth one-periodic map $\widetilde z_i: S^1 = \R/\Z
\to M$ that satisfies
$$
\dot {\widetilde z}_i(t) = X_{H_i'}(\widetilde z_i(t)). \tag 6.7
$$
From the hypothesis $\|H_i - G\|_{C^1} \to 0$ and $T_i \to 1$, we
have
$$
\max_{(t,x) \in [0,1] \times M} |X_{G}(t,x) - X_{H_i'}(t,x)| \to
0. \tag 6.8
$$
It follows from (6.7) and (6.8) that $|\dot z_i|_{C^0} < C$, $C$
independent of $i$. Therefore the sequence $\widetilde z_i:S^1 \to
M$ is equi-continuous and so uniformly converges to a continuous
map $z_\infty: S^1 \to M$ which can be easily shown to be a weak
solution of the equation $\dot z = X_G(z)$. Using the smoothness
of $G$, $z_\infty$ is a genuine one-periodic solution which is
smooth. Therefore by the hypothesis in the Minimality Conjecture,
$z_\infty$ must be a constant solution. This in turn implies that
the sequence $\widetilde z_i$ and so $z_i$ uniformly converges to
a constant map. This already rules out the possibility (i) and
(ii).

For the case of (iii), since $\widetilde{z}_i$ is smooth (and so
$C^1$) and
$$
\text{diam}(\widetilde z_i) = \text{diam}(z_i) \leq
\frac{1}{2}\text{inj}(g)
$$
we can represent $\widetilde  z_i$ in its center of mass
representation
$$
\widetilde z_i = \exp_{x}(\xi_i'), \quad x = x_{\widetilde z_i} =
x_{z_i}.
$$
It follows from the general properties of the exponential map that
we have
$$
Area(w_i') \leq 2 \pi C \max_{t\in [0,1]}|\xi_i'|
$$
with a constant $C$ depending only on $M$. This converges to zero
since $\widetilde z_i$ uniformly converges to a constant map and
so  $\xi_i' \to 0$. On the other hand, we have
$$
\Big|\int_{w_i'} \omega\Big| \leq Area(w_i') \to 0 \tag 6.9
$$
from (6.3). However (6.9) then contradicts to the third
possibility (iii). This finishes the proof. \qed\enddemo

Here we would like to emphasize that in this proof we have used
only $C^1$-small condition for the perturbation, which only gives
rise to the $C^0$-closeness of the corresponding Hamiltonian
vector fields. Recall that a $C^0$-small perturbation of a vector
field does not imply the uniform convergence of the whole flow.
The upshot is that we were only interested in periodic orbits with
periods $1 -\delta < T < 1 +\delta$ which enters in the above
proof in a crucial way.

By this proposition, we have only to achieve the
$\mu$-undertwistedness only for those critical points
$$
[z,w_z] \in \text{\rm Crit}\AA_H
$$
that bifurcate from the critical points of the autonomous $G$. The
main point of the $C^1$ Perturbation Conjecture is that {\it we
hope to be able to change the Conley-Zehnder index arbitrarily
using a $C^1$ small but $C^2$ big perturbation so that the
Conley-Zehnder index for these critical points lie in the region
$-n \leq \mu_H([z,w_z]) \leq n$}. As we mentioned before, the
requirement about $\mu$-undertwistedness would not be possible if
a $C^2$-small perturbation were asked for. This is because both
nondegeneracy of periodic orbits and the Conley-Zehnder index are
stable under $C^2$-small perturbations.

Since all other requirements will be satisfied by any $C^1$ small
perturbation, once we achieve this index condition, we can make
another $C^2$-small perturbation afterwards to achieve the Floer
theoretic nondegeneracy keeping all other properties intact.
Theorem $A$ will then imply that $H$ is length minimizing. We can
apply Theorem 2.3 to prove that $G$ is length minimizing, by
letting the $C^1$ perturbation go to zero. This will finish proof
of the Minimality Conjecture for the very strongly semi-positive
case, once the above $C^1$ Perturbation Conjecture is proven.

One final additional remark is that our conjecture does not
contradict to the Arnold conjecture
$$
HF^*(H) \cong H^*(M;\Lambda_\omega).
$$
It suggests that near the autonomous Hamiltonian as in the
Minimality Conjecture, one can find a $C^1$-close nondegenerate
Hamiltonian whose Floer complex is very special in that its
homologically essential part of the Floer complex resembles that
of the Morse complex of nearby Morse functions.

\head{\bf Appendix: Proof of the index formula (1.5)}
\endhead

As we mentioned before, different sign conventions have been used
in the definitions of various objects in the literature of
symplectic geometry. The only thing that enters in the definition
of the Maslov index is, however, a periodic solution of the
Hamilton's equation
$$
\dot x = X_H(x)
$$
on a symplectic manifold $(M,\omega)$ for a one-periodic
Hamiltonian function $H: S^1 \times M \to \R$. Our convention is that
$X_H$ is defined by
$$
X_H\, \rfloor \omega = dH \quad\text{or equivalently } \,
dH(x)(\xi)=\omega(X_H, \xi). \tag A.1
$$
Furthermore the canonical symplectic form of on $T^*\R^n = \R^{2n}
\cong \C^n$ in the coordinates $z_j = q_j + i p_j$ is given by
$$
\omega_0 = \sum_{j=1} dq^j \wedge dp_j. \tag A.2
$$
This means that on $\R^{2n}$, $X_H = J_0 \nabla H$ where $J_0$ is
the standard complex structure on $\R^{2n} \cong \C^n$ obtained by
multiplication by the complex number $i$. With these being
mentioned, we give the proof of the index formula (1.5) in several
steps.

\medskip

0. There is another package of conventions that have been
consistently used by Polterovich [Po] and others. In that
convention, there are two things to watch out in relation to the
index formula, when compared to our convention. The first thing is
that their definition of the Hamiltonian vector field, also called
as the symplectic gradient and denoted by $\text{sgrad}\, H$, is
given by
$$
\text{sgrad}\, H \, \rfloor \omega = - dH. \tag A.3
$$
Therefore we have $X_H = - \text{sgrad}H$. The second thing is
that their definition of the canonical symplectic form on $T^*\R^n
= \R^{2n} \cong \C^n$ in the coordinates $z_j = q_j + i p_j$ is
given by
$$
\omega'_0 = \sum_{j=1} dp_j \wedge dq^j = -\omega_0\tag A.4
$$
Cancelling out two negatives, the definition of the Hamiltonian
vector field of a function $H$ on $\R^{2n}$ in this package
becomes the same vector field as ours that is given by
$$
J_0 \nabla H
$$
where $\nabla H$ is the usual gradient vector field of $H$ with
respect to the standard Euclidean inner product on $\R^{2n}$.
\smallskip

1. We follow the definition from [CZ], [SZ] of the Conley-Zehnder index
for a paths $\alpha$ lying in $\SS P^*(1)$ where we denote
$$
\SS P^*(1) = \{ \alpha : [0,1] \to Sp(2n,\R) \mid \alpha(0) = id,
\, \det(\alpha(1) - id) \neq 0 \} \tag A.5
$$
following the notation from [SZ]. We denote by $\mu_{CZ}(\alpha)$
the Conley-Zehnder index of $\alpha$ given in [SZ]. Note that the
definition of $Sp(2n,\R)$ are the same in both of the above two
conventions and so the Conley-Zehnder index function $\mu_{CZ}:
\SS P^*(1) \to \Z$ is the same under the above two conventions.
\smallskip

2. A given pair $[\gamma,w] \in \widetilde \Omega_0(M)$ determines
a preferred homotopy class of trivialization of the symplectic
vector bundle $\gamma^*TM$ on $S^1 = \part D^2$ that extends to  a
trivialization
$$
\Phi_w: w^*TM \to D^2 \times (\R^{2n},\omega_0)
$$
over $D^2$ of where $D^2 \subset \C$ is the unit disc with the
standard orientation.
\smallskip

3. Let $z: \R/\Z \times M$ be a one-periodic solution 
of $\dot x = X_H(x)$. Any such one-periodic solution has the form 
$z(t) = \phi_H^t(p)$ for a fixed point $p = z(0) \in \text{Fix}(\phi_H^1)$.
When we are given a one-periodic solution
$z$ and its bounding disc $w:D^2 \to M$, we
consider the one-parameter family of the symplectic maps
$$
d\phi_H^t(z(0)): T_{z(0)}M \to T_{z(t)}M
$$
and define a map $\alpha_{[z,w]}: [0,1] \to Sp(2n,\R)$ by
$$
\alpha_{[z,w]}(t) = \Phi_w(z(t))\circ d\phi_H^t(z(0)) \circ
\Phi_{w}(z(0))^{-1}. \tag A.6
$$
Obviously we have $\alpha_{[z,w]}(0) = id$, and nondegeneracy of
$H$ implies that
$$
\det(\alpha_{[z,w]}(1) -id) \neq 0
$$
and hence
$$
\alpha_{[z,w]} \in \SS P^*(1). \tag A.7
$$
Then the Conley-Zehnder index of $[z,w]$ is, by definition, given
by
$$
\mu_H([z,w]) : = \mu_{CZ}(\alpha_{[z,w]}). \tag A.8
$$
\smallskip

4. When we are given two maps
$$
w, \, w': D^2 \to M
$$
with $w|_{\part D^2} = w'|_{\part D^2}$, we define the glued map
$u = w\# \overline w': S^2 \to M$ in the following way:
$$
u(z) = \cases w(z) & \quad z \in D^+ \\
w'(1/\overline z) & \quad z \in D^-.
\endcases
$$
Here $D^+$ is $D^2$ with the same orientation, and $D^-$ with the
opposite orientation. This is a priori only continuous but we can
deform to a smooth one without changing its homotopy class by
`flattening' the maps near the boundary: In other words, we may
assume
$$
w(z) = w(z/|z|) \quad \text{for } |z| \geq 1 -\e
$$
for sufficiently small $\e >0$. We will always assume that the
bounding disc will be assumed to be flat in this sense. With this
adjustment, $u$ defines a smooth map from $S^2$.
\smallskip

5. For the given $[z,w], \, [z,w']$ with a periodic solution $z(t)
= \phi_H^t(z(0))$, we impose the additional {\it marking}
condition
$$
\Phi_w(z(0)) = \Phi_{w'}(z(0)) \tag A.9
$$
as a map from $T_{z(0)}M$ to $\R^{2n}$ for the trivializations
$$
\Phi_w, \Phi_{w'} : w^*TM \to D^2 \times (\R^{2n},\omega_0)
$$
which is always possible. {\it With this additional condition}, we
can write
$$
\alpha_{[z,w]}(t) = S_{ww'}(t)\cdot \alpha_{[z,w']}(t) \tag A.10
$$
where $S_{ww'}: S^1 =  \R/\Z \to Sp(2n,\R)$ is the {\it loop}
defined by the relation (A.10). Note that this really defines a
loop because we have
$$
\align \alpha_{[z,w]}(0) & = \alpha_{[z,w']}(0) (= id) \tag A.12\\
\alpha_{[z,w]}(1) & = \alpha_{[z,w']}(1) \tag A.13
\endalign
$$
where $(A.13)$ follows from the marking condition (A.9). In fact,
it follows from the definition of (A.5) and (A.9) that we have the
identity
$$
\align S_{ww'}(t) & = \Big(\Phi_w(z(t))\circ d\phi_H^t(z(0))\circ
\Phi_w(z(0))^{-1}\Big) \\
& \qquad \circ \Big(\Phi_{w'}(z(t))\circ
d\phi_H^t(z(0))\circ \Phi_{w'}(z(0))^{-1}\Big)^{-1} \\
& = \Phi_w(z(t))\circ \\
& \qquad \Big(d\phi_H^t(z(0))\circ \Phi_w(z(0))^{-1}\circ
\Phi_{w'}(z(0))\circ(d\phi_H^t)^{-1}(z(0))\Big) \\
& \qquad \circ (\Phi_{w'}(z(t)))^{-1}. \tag A.14
\endalign
$$
Then the marking condition (A.9) implies the middle terms in
(A.14) are cancelled away and hence we have proved
$$
S_{ww'}(t) =\Phi_w(z(t))\circ \Phi_{w'}(z(t))^{-1} \tag A.15
$$
Then we derive the following formula, from the definition
$\mu_{CZ}$ in [CZ] and from (A.10),
$$
\mu_{CZ}(\alpha_{[z,w]}) = 2 \text{ wind}(\widehat S_{ww'})
+\mu_{CZ}(\alpha_{[z,w']}) \tag A.16
$$
where $\widehat S_{ww'}: S^1 \to U(n)$ is a loop in $U(n)$ that is
homotopic to $S_{ww'}$ inside $Sp(2n,\R)$. Such a homotopy always
exists and is unique upto homotopy because $U(n)$ is a deformation
retract to $Sp(2n,\R)$. And $\text{wind}(\widehat S_{ww'})$ is the
degree of the obvious determinant map
$$
{\det}_\C(\widehat S_{ww'}): S^1 \to S^1.
$$
\smallskip

6. Finally, we recall the definition of the first Chern class
$c_1$ of the symplectic vector bundle $E \to S^2$. We normalize
the Chern class so that the tangent bundle of $S^2$ has the first
Chern number 2, which also coincides with the standard convention
in the literature. We decompose $S^2 = D^+ \cup D^-$, and consider
the trivializations $\Phi_+: E|_{D^+} \to D^2 \times
(\R^{2n},\omega_0)$ and $\Phi_-: E|_{D^-} \to D^2 \times
(\R^{2n},\omega_0)$. Denote by the transition matrix loop
$$
\phi_{+-}: S^1 \to Sp(2n,\R)
$$
which is the loop determined by the equation
$$
\Phi_+|_{S^1}\circ (\Phi_-|_{S^1})^{-1}(t,\xi) = (t,
\phi_{+-}(t)\xi)
$$
for $(t,\xi) \in E|_{S^1}$, where $S^1 = \part D^+ = \part D^-$.
Then, by definition, we have
$$
c_1(E) = \text{wind}(\widehat \phi_{+-}) \tag A.17
$$
Now we apply this to $u^*(TM)$ where $u=w\# \overline w'$ and
$\Phi_w$ and $\Phi_{w'}$ are the trivializations given in 4. It
follows from (A.15) that $S_{ww'}$ is the transition matrix loop
between $\Phi_w$ and $\Phi_{w'}$. Then by definition, the first
Chern number $c_1(u^*TM)$ is provided by the number
$\text{wind}(\widehat S_{ww'})$ of the loop of unitary matrices
$$
\widehat S_{ww'}: t\mapsto \widehat S_{ww'}(t); \quad S^1 \to
U(n). \tag A.18
$$
One can easily check that this winding number, {\it not that of
the inverse loop $\widehat S_{ww'}^{-1}$}, is indeed 2 when
applied to the tangent bundle of $S^2$ and so consistent with the
convention of the Chern class that we are adopting.
\smallskip

7. Combining these steps, we have finally proved

\proclaim{Theorem C} Let $(M,\omega)$ be a symplectic manifold and
$X_H$ a Hamiltonian vector field defined by
$$
X_H \rfloor \omega = dH
$$
of any one-periodic Hamiltonian function $H: [0,1] \times M \to
\R$. For a given one-periodic solution $z: S^1 = \R/\Z \to M$ of
$\dot x = X_H(x)$ and two given bounding discs $w, \, w'$, we have
the identity
$$
\mu_H([z,w]) = \mu_{H}([z,w']) + 2 c_1([w\# \overline w']).
$$
The same formula holds for the other package of conventions given
in the paragraph 0 without change.
\endproclaim

\head {\bf References}
\endhead

\enddocument